\documentclass[11pt]{article}
\usepackage[cp1251]{inputenc}
\usepackage[T2A]{fontenc}
\usepackage[english]{babel}
\usepackage{amsfonts}
\usepackage{amsmath}
\usepackage{amssymb}
\usepackage[dvips]{graphicx}

\newcommand\R{{\mathbb R}}
\newcommand\T{{\mathbb T}}
\newcommand\Z{{\mathbb Z}}
\newcommand\Q{{\mathbb Q}}
\newtheorem{theorem}{Theorem}
\newtheorem{lemma}{Lemma}
\newtheorem{proposition}{Proposition}

\newtheorem{definition}{Definition}

\newtheorem{hypothesis}{Hypothesis}

\title{Geometry of symplectic partially hyperbolic \\automorphisms on 4-torus}
\author{L.M. Lerman, K.N. Trifonov\\
\normalsize
Dept. of Differ. Equations, Calculus and Numerical Analysis\\
\normalsize
Lobachevsky State University of Nizhny Novgorod}
\date{}
\begin{document}
\maketitle

\begin{abstract}
We study topological properties of automorphisms of 4-dimensional
torus generated by integer symplectic matrices. The main
classifying element is the structure of the topology of a foliation generated by
unstable leaves of the automorphism. There two different case, transitive
and decomposable ones. For both case the classification is given.
\end{abstract}

\section{Introduction}

We study topological properties of automorphisms of 4-dimensional
torus generated by integer symplectic transformations of $\R^4$. Usually such transformations are
called symplectic automorphisms of the torus. Our main concern is about
partially hyperbolic case. There are many results about partially
hyperbolic diffeomorphisms \cite{Ham_Pot,Hertz,Burns_Wilk,HassPesin} since the basic paper
\cite{Brin_Pesin}. As for partially hyperbolic automorphisms on a torus $\T^n$ is concerned,
the most detailed their study was done in \cite{Hertz} solving the question on their stable
ergodicity posed in \cite{HPS}. Recall that stable ergodicity of a $C^r$-smooth
diffeomorphism $f$ of a manifold $M$ ergodic with respect to a smooth Lebesque measure on $M$
means the existence of neighborhood $\mathcal U\ni f$ in the space of $C^r$-smooth diffeomorphisms
such that any $g\in \mathcal U$ is ergodic. of Our goal here is to classify possible types of the
orbit behavior of symplectic automorphisms on $\mathbb T^4$.
We hope this will be useful as good examples.

Consider standard torus $\mathbb T^n = \mathbb R^n/\Z^n $ as the factor
group of the abelian group $\R^n$ with respect to its discrete subgroup $\Z^n$ of integer
vectors. Denote $p: \R^n \to \T^n$ the related group homomorphism being
simultaneously a smooth covering map.
The coordinates in the space $\R^n$ will be denoted by $x=(x_1,\ldots,x_n).$ Let $A$ be
an unimodular matrix with integer entries. Since the linear map $L_A:x\to Ax$, generated by the matrix
$A$, transforms the subgroup $\Z^n$ onto itself, such a matrix generates diffeomorphism $f_A$ of the
torus $\mathbb{T}^n$ called the automorphism of the torus \cite{Anosov,AS,Franks}. Topological
properties of such maps are the classical object of research (see, for example, \cite{Anosov,
Franks,KH,Manning}). Because the torus automorphism also preserves the standard volume element
 $dx_1\wedge\cdots \wedge dx_n$ on the torus carried over from $\R^n$, then its ergodic properties
have also been the subject of research in many works \cite{Halmos,Bowen, Sinai}. The following
classical Halmos theorem holds for automorphisms of a torus \cite{Halmos}.
\begin{theorem}[Halmos]
A continuous automorphism $f$ of a compact abelian 	group $G$ is ergodic (and mixing) if and only
if the induced automorphism on the character group $G^*$ has not finite orbits.
\end{theorem}

Recall \cite{Pont} the group of characters $G^*$ of a topological group $G$ be the set of
its homomorphisms $ \chi:G\to S^1$ into the group of rotations of the circle $ S^1=\R/\Z$ endowed
with the product of two such homomorphisms $\chi_1,\chi_2$ as follows: $(\chi_1*\chi_2)(g)=
\chi_1(g)+\chi_2(g)\; (\mbox{\rm mod}\; 1).$ If $f:G\to G$ is a group automorphism, then the induced
automorphism $ f^*: G^*\to G^*$ on the character group is defined as $[f^*\chi](g)=\chi(f(g)).$
If a topological group is compact, its character group is discrete \cite{Pont}.

In the case of the abelian group $\mathbb{T}^n=\R^n/\Z^n$, the characters are complex-valued
functions $\chi_m(x)=\exp[2\pi i(m,x)]$, where $m=(m_1,\ldots,m_n)\in \Z^n$ and $(m,x)$ is the
standard inner product. If $L_A$ is the automorphism generated by the matrix $A$, then the induced
automorphism on the character group acts as follows
$$
[L^*_A(\chi_m)](x)= \exp[2\pi i (m,Ax)]= \exp[2\pi i (A^\top m,x)],
$$
where $A^\top$ means the transposed matrix. The existence of
a finite orbit of the induced automorphism means that for some vector $m_0\in\Z^n$ and
natural $k$ the identity holds
$$
\exp[2\pi i ((A^k)^\top m_0,x)]\equiv \exp[2\pi i (m_0,x)]
$$
for all $x\in \mathbb T^n$. This means the integer vector $m_0$ be an eigen-vector of matrix
$(A^\top)^k$ with the eigenvalue 1.
Thus, among eigenvalues of $A^\top$ there exists $\lambda$ such that $\lambda^k =1$. But eigenvalues
of $A$ and $A^\top$  are the same. Conversely, if an integral unimodular matrix $A$ has
an eigenvalue $\lambda$ being a root of unity, $\lambda^k =1,$ then there exists a nonzero vector
$x$ such that $(A^k - E)x = 0$. But matrix $A^k - E$ has integer entries and all its minors
be integer numbers.
Therefore, vector $x$ can be chosen with rational entries and hence can be made an integer vector.
Thus, the transpose matrix $A^\top$ has eigenvalue $\lambda$ with an
integer eigenvector $m_0,$ $(A^\top)^k m_0 = m_0$, and the related
induced action on the character group has a finite orbit. Finally,
by the Halmos theorem, automorphism $L_A$ is ergodic if and
only if the matrix $A$ has not eigenvalues being roots of unity.$\blacksquare$

In fact, the following assertions are equivalent \cite{KK}:

(1) the automorphism $f_A$ is ergodic with respect to Lebesgue measure;

(2) the set of periodic points of $f_A$ coincides with the set of points in $\mathbb T^N$
with rational coordinates;

(3) none of the eigenvalues of the matrix $A$ is a root of unity;

(4) the matrix $A$ has at least one eigenvalue of absolute value greater than one and has
not eigenvectors with rational coordinates;

(5) all orbits of the dual map $A^*: \Z^N \to \Z^N$, other than the trivial zero orbit,
are infinite.

One more result about torus diffeomorphisms is due to Bowen \cite{Bowen} and
allows one to calculate the topological entropy of such automorphism.
\begin{theorem}(Bowen)
If $L_A:\mathbb{R}^n\to\mathbb{R}^n$ is a linear map generated by unimodular integer matrix $A$,
then $$h_d(f_A)=\sum_{|\lambda_i|>1}
\log|\lambda_i|,$$ where $\lambda_1,\dots, \lambda_n$ are the eigenvalues of $A$.
\end{theorem}
In particular, for a partially hyperbolic automorphism of $\T^4$
(see below) the topological entropy is always positive and equals
to $\log|\lambda|,$ $|\lambda| > 1.$
So such automorphism always possesses some type of chaoticity.

If the dimension of the torus is even, $n=2m$, then one can introduce the standard symplectic
structure on $\T^{2m}$ using coordinates in $\R^{2m}$: $\Omega = dx_1\wedge dx_{m+1} + \cdots +
dx_m \wedge dx_{2m}$ and consider symplectic automorphisms of the torus which preserve this
symplectic structure. A symplectic automorphism $f_A$ is then defined by a symplectic matrix $A$
with integer entries. Such matrices satisfy the identity $A^\top I A=I$, where a skew-symmetric
matrix $I$ has the form
$$
I = \begin{pmatrix}0&E\\-E&0 \end{pmatrix}.
$$
The identity above implies the product of two symplectic matrices and the inverse matrix of
a symplectic matrix be symplectic, i.e. symplectic matrices form a group denoted by
$\mbox{\rm Sp}(2n,\R)$ w.r.t. the operation of matrix multiplication. This group is one of the
standard matrix Lie groups \cite{Chevallier}.

Recall the well-known statement about the characteristic polynomial of a symplectic matrix
(see, for instance, \cite{Arn}).
\begin{proposition}
Characteristic polynomial of a symplectic matrix is reciprocal $\chi(\lambda)=
\lambda^{2n}\chi(1/\lambda).$ If $\lambda$ is an eigenvalue of a real symplectic matrix $A$,
then the numbers $\lambda^{-1}, \bar{\lambda}$, $\bar{\lambda}^{-1}$ are also its eigenvalues,
all they have the same multiplicity and the same structure of elementary divisors.
Eigenvalues of $\lambda=\pm 1$ have an even multiplicity, their elementary divisors
of an odd order, if any, and meet in pairs.
\end{proposition}

Another (non-standard) symplectic structure on the torus $\T^{2m}$ can also be defined.
To this purpose, let us choose a non-degenerate skew-symmetric $2m\times 2m$ matrix $J$,
$J^\top = - J$. Such matrix defines a bilinear 2-form
$[x,y]=(Jx,y)$ in $\R^{2m}$, where $(\cdot,\cdot)$ is the standard coordinate inner product.
This form is also called the skew inner product \cite{Arn}. Then a linear map $S:\R^{2m}\to\R^{2m}$
is called symplectic, if the identity $[Sx,Sy]=[x,y]$ holds for any $x,y\in\R^{2m}$. Using
the representation of skew inner product via the matrix $J$ and properties of the
inner product, we obtain the following identity for the matrix $S$ of the symplectic map:
$S^\top JS=J$. This construction allows one to define a symplectic structure on the torus as well,
if matrices $J,S$ are unimodular integer ones. The unimodularity of $S$ follows from its
symplecticity.
Then on the torus a symplectic 2-form is given and map $S$ defines a symplectic automorphism
with respect to this symplectic form. For example, let $B$ be any nondegenerate integer unimodular
matrix. Having the standard skew inner product $(Ix,y)$ in $\R^{2m}$ we can define new skew
inner product as $[x,y]= (IBx,By) = (B^\top I Bx,y)$. Since the matrix $J = B^\top I B$ is skew
symmetric, integer and nondegenerate, this skew inner product generates symplectic 2-form on the torus.

Let  $P(x)=\lambda^4+a\lambda^3+b\lambda^2+a\lambda+1$ be some reciprocal polynomial with integer
coefficients which is irreducible over the field $\mathbb{Q}$ and suppose it to have two complex
conjugate roots on the unit circle and two real eigenvalues outside of the unit circle. The companion
matrix of this polynomial $Q$ has the form
$$
A=\begin{pmatrix}0&1&0&0\\0&0&1&0\\0&0&0&1\\-1&-a&-b&-a\end{pmatrix}.
$$
This matrix is not symplectic $A^\top I A \ne I$ with respect to the standard symplectic
2-form $[x,y]=(Ix,y)$, $x,y,\in \R^4$. Let us show, however, $A$ to be symplectic with respect to
a nonstandard symplectic structure on $\mathbb{R}^4$ defined as $[x,y]=(Jx,y),$ $x,y\in \mathbb {R}^4$
with a skew-symmetric integer non-degenerate matrix $J$.
We consider the identity $A^\top JA=J$ as the algebraic system $A^\top J-JA^{-1}=0$ for
the entries of the matrix $J$. A solution to this matrix equation gives
$$
J=\begin{pmatrix}0&0&1&0\\0&0&-a&1\\-1&a&0&0\\0&-1&0&0\end{pmatrix}.
$$
Thus, the matrix $A$ becomes symplectic with respect to this nonstandard symplectic
structure. This will be used later on.

The structure of the paper is as follows. In Section 2 we present possible types symplectic
automorphisms on a symplectic torus when matrices $A$ have not unit eigenvalues. In particular,
we distinguish among them partially hyperbolic automorphisms. In Section 3 we find invariant
foliations related to such automorphisms. In Section 4 we show the existence of partially hyperbolic
symplectic automorphisms with transitive unstable foliations. In section 5 we obtain
a classification of symplectic toral automorphisms.

\section{Symplectic automorphisms on $\mathbb T^4$}

Consider the standard four-dimensional torus $\mathbb{R}^4/\mathbb{Z}^4=\mathbb{T}^4$. We assume
a matrix $A$ to be symplectic, $A\in Sp(4,\mathbb{Z}).$ Then the automorphism $f_A$ of
the torus is symplectic with respect to standard symplectic 2-form $dx_1\wedge dy_1+dx_2\wedge
dy_2$ w.r.t. the standard coordinates $(x_1,y_1,x_2,y_2)$ in $\R^4.$

Consider two general classes of real symplectic matrices.
\begin{itemize}
\item  The symplectic integer matrix $A$ is hyperbolic, i.e., $A$ has not eigenvalues on the unit
circle. Therefore, its eigenvalues form two pairs  $\lambda_1, \lambda_1^{-1}$, $\lambda_2,
\lambda_2^{-1}$, and $|\lambda_i| \ne 1.$ Here two cases are possible: either both pairs are
real and $|\lambda_i|<1$
(the saddle case), or they form a complex quadruple $\rho\exp[\pm i\alpha]$,
$\rho^{-1}\exp[\pm i\alpha]$, $0<\rho<1,$ (the saddle-focus case). In any case, the diffeomorphism
on the torus is an Anosov automorphism of the torus \cite{Anosov,AS,KH}. Matrix $A$ defines in
$\mathbb{R}^4$ a linear map $L_A$ with a fixed point at the origin $O$. This point is a saddle or
saddle-focus, respectively. For such fixed point there exist two-dimensional stable and unstable
manifolds, here they are two-dimensional subspaces in $\mathbb{R}^4$ being invariant under
the action of $L_A$, they intersect each other transversally at $O$.
Stable subspace $W^s$ consists of vectors in $\mathbb{R}^4$ which are contracted exponentially
under iterations $L_A^n$ as $n\to\infty$.
Similarly, unstable subspace $W^u$ consists of vectors in $\mathbb{R}^4$ which are
contracted exponentially under iterations $L_A^n$ as $n\to -\infty$.

When projecting on the torus, the subspace $W^s$ and all affine 2-dimensional
planes, being factor classes of $\R^4/W^s$, form a foliation of the torus whose
leaves are embedded two-dimensional planes. The same is true for unstable subspace $W^u$: when
projecting onto a torus, it generates the unstable foliation on the torus. Each leaf of stable
foliation is dense on a torus and intersects transversally each leaf of the unstable foliation, whose
leaves are also dense on the torus.

\item The symplectic map $L_A$ in $\R^4$ is partially hyperbolic that corresponds to matrices $A$
with a pair of eigenvalues on the unit circle $\exp[\pm i\alpha]$ and a pair of real numbers
$\lambda, \lambda^{-1},$ $|\lambda|<1.$ Such linear map has a fixed point at the origin in
$\mathbb{R}^4$ called 1-elliptic point \cite{LMar}. This fixed point has a one-dimensional stable
eigenspace $l^s$ corresponding to eigenvalue $|\lambda|<1$, one-dimensional unstable eigenspace $l^u$
corresponding to the eigenvalue $|\lambda^{-1}|>1$, two-dimensional invariant center subspace
$W^c$ corresponding to the eigenvalues on the unit circle. The projection $p:\R^4\to \T^4$
allows, using shifts on the torus, to transfer these subspaces into the tangent spaces
at any point of the torus. Then  the tangent space $T_x\T^4$ of the torus at each point $x$ splits
into a direct sum of three subspaces (two in the hyperbolic case, stable and unstable ones)
$$
T_x{\mathbb T^4} = E_x^s \oplus E_x^c \oplus E_x^u,
$$
where for vectors from $E_x^s$ there is an exponential contraction under iterations of the differential
$DL_A$ (in fact, the matrix $A$), on $E_x^u$ there is an exponential expanding,
and vectors in $E_x^c$ are uniformly bounded under iterations with respect to any Riemannian
metric on the torus.
\end{itemize}

The characteristic polynomial of a symplectic matrix has the form $\lambda^4+a\lambda^3+b\lambda^2+
a\lambda+1.$ Therefore, if we denote $\mu=\lambda+\lambda^{-1}$, then
the characteristic polynomial $\chi(\lambda)$ is expressed as $\chi(\lambda)= \lambda^2 (\mu^2 +
a\mu + b-2).$ In the case when symplectic matrix $A$ generates a partially hyperbolic mapping,
the value of $\mu_1 = \lambda + 1/\lambda$, $0<|\lambda|<1,$ is either greater than 2 or
lesser than -2, and the value $\mu_2 = \exp[i\alpha]+\exp[-i\alpha] = 2\cos\alpha$ does not
exceed 2 in modulus. So, in this case one has $|\mu_1| > 2,$ $|\mu_2| < 2,$ if
$\alpha \ne 0,\pi\;(\rm\;mod\;2\pi)$.

\section{One-dimensional foliations \\of partially hyperbolic automorphisms}

Now we consider partially hyperbolic automorphisms of the torus. Partially hyperbolic
diffeomorphisms on any smooth manifold $M$ were first introduced and studied in
\cite{BP}. Here we use a modification of the definition in \cite{Hertz}. Let $L$ be a linear
transformation between two normed linear spaces. The norm,
respectively conorm, of $L$ are defined as
$$||L||:= \sup{||Lv||, ||v|| = 1},\; m(L):= \inf{||Lv||, ||v|| = 1}.$$
\begin{definition}\label{parhyp}
A diffeomorphism $f:M\to M$  is partially hyperbolic, if there
is a continuous $Df$-invariant splitting
$TM=E^u\oplus E^c\oplus E^s$ in which $E^u$ and $E^s$ are nontrivial sub-bundles and
  $$ m(D^uf)>||D^cf||\ge m(D^cf)>||D^sf||,$$
$$m(D^uf)>1>||D^sf||,$$ where $D^\sigma f$ is the restriction of
$Df$ to $E^\sigma$ for $\sigma=s,c$ or $u$.
\end{definition}
In the case of a 4-dimensional torus, a symplectic partially hyperbolic automorphism is defined by
a symplectic matrix having two (simple) real eigenvalues
$\lambda,\lambda^{-1}$ outside the unit circle
and two complex conjugate eigenvalues on the unit circle. The corresponding eigenspaces
$W^s, W^u$ and the subspace $W^c$ when projecting on a torus and shifting them at any point
give the required decomposition of the definition.

When studying the geometry of partially hyperbolic maps with different dynamics, it is useful
to study first the possible behavior of projections onto a torus of eigenspaces for real
eigenvalues $\lambda,\lambda^{-1}.$ In the space $\mathbb{R}^4$ we have the following orbit structure
for the linear map $L_A$ generated by the matrix $A$. Recall that in the case of a symplectic linear
map with a 1-elliptic fixed point $O$, there is a two-dimensional center invariant subspace $W^c$
corresponding to a pair of eigenvalues $\nu, \bar{\nu}, \nu =\exp [i\alpha]$ on the unit circle.
If $\alpha/2\pi\ne r,$ $r\in\mathbb{Q}$ (the non-resonance case, this is equivalent to $\nu^n \ne 1$
for any $n\in \Z$), then the subspace $W^c$ is foliated into closed
invariant curves. This follows from the fact that the restriction of the map $L_A$ on $W^c$ has
a quadratic positive definite invariant function (integral) whose level lines are closed invariant
curves (in the resonance case $\alpha/2\pi = r$ a quadratic integral also exists but all its level
lines consist
of periodic points of the same period, due to linearity, therefore one can construct any invariant
sets from them). The restriction of the map $L_A$ to any such curve is conjugated to a rotation
on the circle $\varphi\to\varphi+ \alpha$ (mod $2\pi$) by the angle $\alpha$, and the rotation number
is independent of the curve due to the linearity of the map. For the case of irrational number
$\alpha/2\pi$ these invariant curves are defined correctly and the shift on each such curve is
transitive.

In addition to the indicated invariant lines and the center plane, there are two more three-dimensional
invariant subspaces in $\R^4$ spanned, respectively, by vectors from subspaces $W^c$ and $l^s$
(the center-stable 3-plane $W^{cs}$) and vectors from $W^c$ and $l^u$ (the center-unstable
3-plane $W^{cu}$). The factor-classes $\mathbb{R}^4/W^{cu}$ and $\mathbb{R}^4/W^{cs}$ define
two invariant 3-foliations in $\R^4$ into three-dimensional affine planes. Here, the invariance is
understood in the following sense: the image with respect to $L_A$ of a leaf of the foliation
coincides with some leaf (possibly with another one) of the same foliation. All orbits of $L_A$
not lying in the union $W^{cs}\cup W^{cu},$ go to infinity for both positive and negative iterations
of $L_A$. In particular, the mapping $L_A$ has not any other three-dimensional invariant
subspaces except those two $W^{cs}$ and $W^{cu}$. If $\exp [i\alpha]$ is not a root of unity, then
subspaces $W^{cs}$, $W^{cu}$ are foliated  into two-dimensional cylinders being stable
(respectively, unstable) invariant manifolds of invariant curves on the central plane $W^c$.

The projection of the eigen-lines $l^s, l^u$ onto the torus can lead to different situations.
To understand  this it is necessary to clarify what is a projection of a subspace from
$\mathbb{R}^4$ on the torus $\mathbb{T}^4=\mathbb{R}^4/\Z^4$. Since the torus is a commutative
Lie group, the tangent space at zero possesses the structure of the commutative Lie algebra,
this is identified with $\mathbb{R}^4$ in the standard way. Then the projection $p$ is an
exponential map of the Lie algebra onto the Lie group. One-dimensional subspace $l^u$ (or $l^s$)
coincides with the one-parameter subgroup $t\gamma^u$ generated by the vector $\gamma^u$,
and its projection is the image under the exponential mapping of the algebra into the group.
This subgroup is included into orbits of the constant vector field on $\mathbb{T}^4$
invariant under shifts. In coordinates $\theta$ on $\mathbb{T}^4$ induced
by coordinates in $\mathbb{R}^4$ we get the vector field $\dot \theta =
\gamma^u.$ Its orbit structure depends on the number of rationally independent
integer solutions of the equation $(\xi,\gamma^u)=0.$ This number can take
values $0,1,2,3$. In the first case, as is known, the related orbits of
the constant vector field are transitive in $\mathbb{T}^4$ (it is a partial case
of the Kronecker theorem, see, for instance, \cite{Levitan}).

The one-parameter subgroup in $\mathbb{T}^4$ generated by $\gamma^u$ is an invariant subset
with respect to the automorphism $f_A$ (this is the strongly unstable curve of the fixed point
$\hat{O}$). Therefore, its closure is also an invariant subset, that is a smooth invariant
torus of some dimension in $\mathbb{T}^4$ \cite{Bourbaki}. As was said above, the dimension
of this torus depends on the number of an integer linear independent relations of the form
$(m,\gamma^u)=0$:
1) none integer relations of the form $(m,\gamma^u)=0,$ $n\in\mathbb{Z}^4$;
2) vector $\gamma^u$ satisfies an only integer
relation (up to multiplication at a constant) $(n,\gamma^u)=0,$ $n\in\mathbb{Z}^4$;
3) there are two such rationally independent relations $(n,\gamma^u)=0,$ $(m,\gamma^u)=0,$
vectors $n,m\in\mathbb{Z}^4,$ are linearly independent; 4) there are three such relations with
linearly independent vectors $m,n,k\in\Z^4$.
In the last three cases, vector $\gamma^u$ is called resonant, in the
first case it is called incommensurate or non-resonant.

In the second case the linear 3-dimensional subspace in $\R^4$ defined by the equation $(n,x)=0$
is projected onto a 3-torus in $\T^4$, and the straight-line spanned by vector $\gamma^u$
does not pass through the points of its integer sub-lattice. Therefore, this line is projected
into a transitive immersed line in this 3-torus.

Similarly, for the third case, the corresponding subspace is
two-dimensional, it is projected onto 2-torus in $\T^4$. The straight-line, spanned by vector
$\gamma^u$, is projected into a transitive immersed line in this 2-torus. In the fourth case,
the straight line necessarily passes through the integer point and therefore it is projected
onto a simple closed curve in $\T^4$ without self-intersections.

First, obviously the eigen-line $l^u$ corresponding to $\lambda > 1$ ($l^s$ for $\lambda < 1$)
can not intersect the integer lattice $\Z^4$. Indeed, if so, the projection in the
torus of this straight line is a closed invariant curve for the map $f_A$ and the restriction
of this map on this circle gives a diffeomorphism on the circle with one unstable (respectively,
stable) fixed point that is impossible. Any trajectory of the vector field
$\dot{\theta}=\gamma^u$ coincides with a leaf of the unstable (stable) invariant foliation
on the torus which consists of projections of all straight lines corresponding to the factor-classes
$\R^4/l^u$ (respectively, $\R^4/l^s$).

The second, the closure of any trajectory of the constant vector field on the torus is
a smooth torus of some dimension, its dimension is equal to 4 minus the number of
rationally independent linear relations for the components of the vector $\gamma^u:$ $(m,\gamma^u)=0$,
$m\in \mathbb Z^4.$ If such relations are absent at all, the vector $\gamma^u$ is incommensurate,
and then each trajectory of the vector field is dense on the torus, due to the Kronecker theorem.

If there is the unique such relation (up to the multiplication at a rational number), the torus
$\mathbb T^4$ is foliated into three-dimensional invariant tori, on each of them any trajectory
of the constant vector field is everywhere dense. If there are two such independent rationally
relations, the torus $\mathbb T^4$ is foliated into two-dimensional tori, on each of them the
trajectory is everywhere dense. The case of three rationally independent relations leads to
the foliation by the circle which in our case is impossible, as was mentioned above.

In fact, for the case we study, the following statement is valid.
\begin{proposition}
The closure of the unstable leaf passing through a fixed point $\hat{O}$ is either
a 2-dimensional torus or the whole $\T^4$, that is, in the latter case the leaf is transitive.
\end{proposition}

{\bf Proof}. Assume the closure of the unstable manifold of a fixed point
$\hat{O}$ on $\mathbb{T}^4$ to form a 3-dimensional torus $\mathbb{T}^3$. This means the
eigenvector $\gamma^u$ to satisfy the unique integer relation
$(m,\gamma^u)=0,$ $m\in\mathbb{Z}^4$. Consider 3-dimensional hyperplane $(m,x)=0$ in $\mathbb{R}^4$
defined by the co-vector $m$. The 3-dimensional torus is the projection of this
hyperplane.

Let us show in this case the stable one-dimensional manifold of the fixed point $\hat{O}$ be
transversal to the torus $\mathbb{T}^3$. The torus $\mathbb{T}^3$ contains the fixed point
$\hat{O}$. Since
the torus is smooth and invariant with respect to the map $f_A$ (as the closure of an invariant set),
the tangent space to the torus at the fixed point is invariant with respect to differential $Df_A$.
The torus $\mathbb{T}^3$ is obtained by the projection of the hyperplane in $\mathbb{R}^4$ passing
through the fixed point $O$ of the map $L_A$. This plane is invariant with respect to $L_A$,
but linear partially hyperbolic map $L_A$ has not other invariant 3-planes through $O$ except for
$W^{cs}, W^{cu}$. Only the second of them
contains the straight line spanned by the vector $\gamma^u.$ So, the torus $\mathbb{T}^3 $ is
a projection of the center-unstable plane. But then the stable eigenvector is transversal to this
plane, therefore the projection onto $\mathbb{T}^4$ of the stable eigen-line is a smooth curve
transversal at the point $\hat{O}$ to the torus $\mathbb{T}^3$. Since
$\mathbb{T}^3$ is a smooth closed submanifold in $\mathbb{T}^4$, there is
a neighborhood $V$ of the point $\hat{O}$ such that all points of the stable
curve $W^s(\hat{O})$ in $V$ do not belong to $\mathbb{T}^3.$

Due to transversality of the torus $\mathbb{T}^3$ and a stable invariant curve at $\hat{O}$,
they must intersect each other at more than one point in $\mathbb{T}^4$ since this curve is not closed.
Stable invariant curve is given as $t\gamma^s\; (\mbox {\rm modd}\;1)$, it contains a point $z$
other than $\hat{O}$ which belongs to $\mathbb{T}^3$. Since $z$ belongs to
$W^{s}$, its backward iterations $f^{-n}_A(z)$ must lie on the stable curve near the point
$\hat{O}$ for positive $n$ large enough. As a consequence, these points do not belong to
$\mathbb{T}^3$ for such iterations. On the other hand, the torus is invariant with respect to $f_A$,
therefore, all iterations of $z$ should lie on it. This contradiction proves that the closure
of an unstable curve is not a 3-torus.

The case of two rationally independent relations for partially hyperbolic matrix is possible. It is
sufficient to choose, for example, a block-diagonal integer matrix $A$ composed of two integer
$2\times 2$-blocks, one of which has eigenvalues $\lambda, \lambda^{-1},$ $|\lambda|< 1$,
and the second block has two complex conjugate eigenvalues $\exp[i\alpha], \exp[-i\alpha]$
on the unit circle. Note in this case, that the characteristic polynomial of such a matrix
is the product of two monic polynomials of second degree
with integer coefficients, i.e. it is reduced over the field of rational numbers and the
numbers $\exp[i\alpha]$ are roots of unity (of degree 3,4,6).$\blacksquare$

For the case when the closure of the unstable curve is a 2-dimensional torus,
the following assertion is valid
\begin{proposition}
If the closure of the unstable (stable) invariant curve of $\hat{O}$ is a two-dimensional torus,
then the characteristic polynomial of $A$ is reducible over rational numbers $\Q$, i.e. it is
the product of two monic polynomials with integer coefficients. In particular, the eigenvalues of
the matrix $A$ lying on the unit circle are roots of unity.
\end{proposition}

{\bf Proof}. Let the closure of an unstable (stable) invariant curve of $\hat{O}$ be a two-dimensional
torus $T$. This torus is a smooth invariant manifold for $f_A$ containing a fixed point
$\hat{O}.$ Tangent plane to $T$ at $\hat{O}$ is the invariant 2-plane w.r.t.
differential $Df_A$. In the covering space $\mathbb{R}^4$ the pre-image through the origin $O$
w.r.t. the projection $p$ of this 2-plane is invariant 2-plane w.r.t. $L_A.$ There are only two such
2-planes: $W^c$ and the plane $W^*$ spanned by two eigen-vectors $\gamma_s, \gamma_u$. Only the
second of them contains $\gamma_u$, so $W^*$ is this plane. There exist two linear independent
integer vectors in $W^*$, since the projection of this plane to
$\mathbb{T}^4$ is the invariant 2-dimensional torus $T$. The lattice generated
by these two integer vectors does not intersect the eigen straight-line, because
its projection is everywhere dense on the torus.
This implies the restriction of $f_A$ to this torus be an Anosov map which is
generated by the restriction of $L_A$ onto invariant 2-plane $W^*$. Eigenvalues of
this restriction are $\lambda, \lambda^{-1}$ with their eigen-vectors $\gamma_s,
\gamma_u.$ In contrast to a transitive case, where in $\mathbb{R}^4$ the
invariant w.r.t. $L_A$ such 2-plane also exists, here this plane contains an integer
sub-lattice with two independent integer vectors.

Consider the minimal polynomial corresponding to the invariant 2-plane $W^*$ in $\R^4$, its degree
is 2. This polynomial can be generated by one integer vector $v\in W^*$, the independent integer
vector $Av$, therefore $A^2v=aAv+bv$, due to the invariance of the plane, $a,b \in \R$. Let us consider
this vector equality as the over-defined system w.r.t. unknowns $a,b.$ The equality says that
the solution of the system exists. Integer vectors $v$, $Av$ are independent, hence there is
a nonzero integer 2-minor in the related coefficient matrix. Vector $A^2v$ is also integer.
Thus existing solution $(a,b)$ of the system consists of rational numbers.

Due to its minimality, the polynomial $z^2 - a z -b$ divides the characteristic polynomial of $A$.
So, the characteristic polynomial of $A$ is the product of two polynomials of degree
two with rational coefficients, i.e.
$$
\lambda^4+\alpha\lambda^3 + \beta\lambda^2+\alpha\lambda + 1 =
(\lambda^2+r\lambda+1)(\lambda^2+s\lambda+1),
$$
with integer $\alpha,\beta$, where $r, s\in \mathbb{Q}$.
Then, equating the coefficients at the same powers, we get the system $r+s=\alpha,\;rs=\beta-2.$
Thus $r$ is a root of the monic quadratic polynomial $t^2 - \alpha t + \beta -2$
with integer coefficients. Hence, its roots are either integer or irrational numbers.
So, $r$ is an integer number and $s$ is as well.
Therefore, the characteristic polynomial splits into the product
of two polynomials with integer coefficients. One of these polynomials,
say $\lambda^2+s\lambda+1=0$ with integer $s$ corresponds to the pair of
roots on the unit circle. This means $|s|<2$ and so $s=-1,0,1$. Then we get
the numbers $\exp[i\alpha]$ be roots of unity of degree 3,4,6.

Now we shall show that the projection onto the torus $\mathbb{T}^4$ of the invariant subspace $W^c$
gives also an invariant two-dimensional torus. To this end, one needs to find two independent
integer vectors in $W^c$. Matrix $A$ has eigenvalues being
roots of unity of degree $k=3,4,6$, $\lambda^k=1$. Then $A^k$ has the double eigenvalue 1.
Therefore, the equation for finding related eigenvectors is of
the form $(A^k-E)x=0$. Matrix $A^k-E$ is integer with all its minors being integer numbers.
Therefore, vector $x$ can be chosen with rational coefficients and hence can be made integer
vector. Thus, matrix $A^k$ has an eigenvalue 1 with an integer eigenvector $m_0,$ $A^k m_0=m_0$.
Notice, that $W^c$ is the invariant plane for $A$, hence $m_0 \in W^c.$
The vector $Am_0$ is also integer and $m_0$ is not an eigenvector, hence $Am_0\ne \gamma m_0$
and vectors $m_0, Am_0$ form the basis of $W^c$. So, the projection of this plane onto $\mathbb{T}^4$
is an invariant two-dimensional torus.
$\blacksquare$

In view of this theorem, we shall call {\em decomposable} the case, when the closure of the
unstable leaf of $\hat{O}$ is a 2-dimensional torus.
Recall the theorem on the rational canonical form \cite{Horn} (Frobenius normal
form).
\begin{theorem}
Every matrix $A$ with rational entries is similar over $\Q$ to a block diagonal matrix
with blocks being companion matrices of its elementary divisors.
\end{theorem}	
It follows that the matrix $A$ of a linear operator $L_A$ is similar to the block matrix $B$
with integer coefficients consisting of two $(2\times 2)$ blocks.
As was shown above, in the decomposable case the characteristic polynomial of the linear operator
$L_A$ splits into the product of two second degree polynomials with integer coefficients.
\begin{proposition}
If the decomposable case is realized, then the matrix $A$ is similar to
$B$ via some integer similarity matrix $T$.
\end{proposition}
{\bf Proof}. By the definition of matrix similarity, the relation $B=TAT^{-1}$ is valid, where
$T$ is the transform matrix. We rewrite this relation in the form $BT-TA=0.$ This matrix equation
in $T$ is equivalent to a system of 16 linear homogeneous equations w.r.t. 16 unknown $T$-entries.
Finding a transforming matrix is reduced to a solution of this system.
Moreover, we have to find among the set of all solutions such that $\det T\ne 0$.
Such a solution exists since the matrices $A$ and $B$ are similar but we need to find an
integer $T$.

The set of solutions of a linear homogeneous system for $T$ is always infinite given
as $T=UT_1$, where $T_1$ is one of transforming matrices, and $U$ is an arbitrary matrix that is
permutable with $B$. This allows one to make the matrix $T$ integer. Since  matrices
$A$ and $B$ are integer, the system $BT-TA=0$ has integer coefficients. As free variables,
we can always choose integers, then the solution is either integer or rational. Thus, the entries
of the matrix $T$ are either integer or rational numbers. In the first case, the statement is proved.
In the second case we use the formula $T=UT_1$, where $T_1$ is a matrix with rational
coefficients, $U$ is a scalar matrix $\mu E$, where $\mu$ the common denominator of matrix
elements for $T_1.$ Matrix $U$ of this kind always permutes with $B$. Then the matrix
$T$ is integer.$\blacksquare$

\section{Partially hyperbolic automorphism of a 4-torus \\ with transitive unstable foliation}

To order to construct examples of symplectic partially hyperbolic automorphisms of a torus
with different dynamical properties we need to find a matrix in $Sp(4,\mathbb Z)$ which have
two complex conjugate eigenvalues on the unit circle and two real eigenvalues
$\lambda$ and $\lambda^{-1}$, $0<|\lambda|< 1$. Moreover, we would like to obtain
an automorphism of the torus whose one-dimensional foliations corresponding to real eigenvalues are
transitive. To obtain such a matrix, we start, following \cite{Kat}, with a irreducible over
rational numbers quadratic polynomial $P(z)$ with integer coefficients having two roots,
one is greater than $2$ and the second is lesser than $2$ in modulus. We make a change of
variable $z=x+x^{-1}$ in this polynomial to obtain a polynomial of the
fourth degree which serves as the characteristic polynomial of the matrix with the desired properties.

For example, we start with the polynomial $P=z^2-3z + 1$ that gives
the polynomial $Q=x^4-3x^3+3x^2-3x+1$ which is irreducible over the field $\mathbb{Q}$.
Then the companion matrix of this polynomial $Q$ \cite{Horn} has the form
$$
A=\begin{pmatrix}0&1&0&0\\0&0&1&0\\0&0&0&1\\-1&3&-3&3\end{pmatrix}.
$$
This matrix possesses the required properties of its eigenvalues. However, this matrix is not
symplectic with respect to the standard symplectic 2-form $[x,y]=(Ix,y)$, $x,y,\in \R^4$:
$A^\top I A \ne I$. We show, however, $A$ to be symplectic with respect to a nonstandard
symplectic structure on $\mathbb{R}^4$ defined as $[x,y]=(Jx,y),$ $x,y\in \mathbb {R}^4$
via a skew-symmetric integer non-degenerate matrix $J$.

To get symplecticity $(JAx,Ax)=(A^\top JAx,x)$, we need the equality $A^\top JA=J$ to hold.
As above, we consider this equality as the algebraic system $A^\top J-JA^{-1}=0$ for the entries of
the matrix $J$. The solution to this matrix equation gives
$$
J=\begin{pmatrix}0&0&1&0\\0&0&-3&1\\-1&3&0&0\\0&-1&0&0\end{pmatrix}.
$$
Thus, the matrix $A$ becomes symplectic with respect to this nonstandard symplectic structure.

Real eigenvalues of matrix $A$ are
$$
\lambda=\frac{3+\sqrt{5}+\sqrt{6\sqrt{5}-2}}{4}> 1,\;\lambda^{-1}.
$$
The eigenvector corresponding to $\lambda$ is $\gamma^u = (1, \lambda, \lambda^2, \lambda^3)^{\top}$.
Factor-classes $\mathbb {R}^4/l^u$ form an invariant foliation under $L_A$ into affine lines.
These lines are projected onto $\mathbb{T}^4$ as trajectories of the vector field
$$
\dot x_1 = 1,\;\dot x_2 = \lambda,\;\dot y_1 = \lambda^2,\;\dot y_2 = \lambda^3.
$$

To prove the orbits of this vector field to be transitive in $\mathbb{T}^4$,
we need to verify the vector $\gamma^u$ to be incommensurate, i.e. $(m,\gamma^u)\ne 0$ for
any integer vector $m\in\mathbb{Z}^4$. The number $\lambda$ is an algebraic number
being the root of the polynomial with integers coefficients.

Recall a number $\xi\in \mathbb C$ be {\em algebraic} if it is a root of a
polynomial with rational coefficients, the number $\xi$ is called {\em integer algebraic},
if it is a root of a polynomial with integer coefficients whose leading coefficient
equals to unity (a monic polynomial). With any algebraic number the notion of its degree
is associated. Namely, if the number $\xi$ is the root of a polynomial of some degree,
then multiplying this polynomial at another polynomial with rational coefficients
we obtain a polynomial of a greater degree with the root $\xi$.
The polynomial of the least degree with rational coefficients having the root $\xi$
is called the {\em minimal} polynomial of the algebraic number $\xi$, and its degree
is called the degree of the number $\xi$. In particular, all rational numbers are algebraic
numbers of degree one. The polynomial $Q$ above has degree 4 due to the
following theorem \cite{Horn}
\begin{theorem}
Every monic polynomial is both the minimal polynomial and the
characteristic polynomial of its companion matrix.
\end{theorem}

\begin{theorem}
Let $\xi$ be an algebraic number over field $\mathbb Q$ and $p$ its minimal monic polynomial.
Then
\begin{enumerate}
\item the polynomial $p$ is irreducible over $\mathbb Q$;
\item the polynomial $p$ is unique;
\item if $\xi$ is a root of some polynomial $f$ over the field $\mathbb Q$ then
$f:p$.
\end{enumerate}
\end{theorem}
From the last statement of this theorem it follows that if the number $\lambda$ is a root of
the irreducible polynomial $Q$ of the fourth order, it could not to be a root of a polynomial
of smaller degree with rational (integer) coefficients.

Now we can prove the lemma
\begin{lemma}\label{den}
The vector $\gamma$ is incommensurate.
\end{lemma}
{\bf Proof.} Suppose vector $\gamma$ be commensurate. Then an integer vector $(m_1,m_2,m_3,m_4)$
exists such that the equality $m_1 + m_2 \lambda + m_3\lambda^2 + m_4\lambda^3
=0$ holds, i.e. $\lambda$ is a root of the polynomial $P$ of the third (or lesser)
degree with integer coefficients. So $P$ divides $Q$ and $Q$
is reducible.$\blacksquare$

Next proposition is well known for a hyperbolic case, its proof for the
partially hyperbolic case follows the lines of \cite{KH}.
\begin{proposition}
The set of periodic points of the map $L_A$ is dense. If unstable the one-dimensional foliation
of the automorphism $f_A$ is transitive, then $f_A$ is a topologically mixing transformation.
\end{proposition}

Before the proof, recall the definition of topologically mixing
map $f$ of the metric space $M$, $f:M\to M$ \cite{KH}.
\begin{definition}
A map $f:M\to M$ is called topologically mixing, if for any two nonempty open sets $U,V\subset M$
there exists a positive integer $N=N(U,V)$ such that for every $n>N$ the intersection $f^n(U)\cap V$
is nonempty.
\end{definition}

{\bf Proof.} The fact that points with rational
coordinates are the only periodic points for $L_A$ and hence the density of periodic orbits
is similar to the proof in the hyperbolic case (see, for instance, \cite{KH}).
Now we prove the topological mixing for map $f_A$, if its unstable foliation is transitive.
Consider any two nonempty open subsets of $U,V\subset\mathbb{T}^4.$ Let $p\in U,q\in V$ be periodic
points and $n$ their common period. Consider 3-dimensional leaf of the center-stable foliation
through the point $p$. This leaf is $f^{n}_A$-invariant. We also consider a 1-dimensional leaf
of the unstable foliation for a point $q$; it is also $f^{n}_A$-invariant.
By the proposition, the leaf of the unstable foliation is everywhere dense in
$\mathbb{T}^4$. In the covering space $\mathbb{R}^4$ subspaces $l^u$ and $W^{cs}$ are
transversal and intersect each other at only point $O$, therefore the 1-dimensional unstable
and 3-dimensional center-stable leaves on the torus are also transversal. Thus there
exists a neighborhood $G$ of the point $p$, in which the unstable foliation generates a smooth
flow box into segments of leaves, and the center-stable foliation generates a smooth
foliation into 3-disks, and leaves of these two foliations are transversal. The intersection
with the neighborhood $G$ of the leaf through $q$ forms an dense set of segments. Segments
of this leaf intersect the three-dimensional disk of the leaf through $p$ at the dense set
of points. Therefore, in the 3-disk there are points of the 1-dimensional leaf accumulating to
$p.$ So, there exists a sufficiently large number $k$ such that $f^{kn}_A(U)\cap V$ is nonempty.
$\blacksquare$

In the example above, the stable 1-dimensional foliation will also be transitive,
because it is generated by the eigenvector $(1,\lambda^{-1},\lambda^{-2},\lambda^{-3})$ and the
integer relation $m_0 + m_1 \lambda^{-1}+m_2\lambda^{-2}+m_3\lambda^{-3}=0$ is impossible.
This fact is not occasional.
\begin{proposition}
If $f_A$ is a symplectic automorphism of $\mathbb{T}^4$ with the dense
unstable foliation, then its stable foliation is also dense.
\end{proposition}

{\bf Proof}. Indeed, let $\gamma_u, \gamma_s$ be related eigen-vectors and
the projection of the line $t\gamma_u$ into $\mathbb{T}^4$ is dense. Then
characteristic polynomial of $A$ is irreducible. Suppose the closure of the
projection of $t\gamma_s$ is a 2-dimensional invariant torus in
$\mathbb{T}^4$. As above, it implies the characteristic polynomial of $A$
is reducible, that is a contradiction.$\blacksquare$

It would be useful to get some general statement concerning the
foliation generated by the center plane $W^c$ for the transitive case. The
projection of $W^c$ onto the torus $\T^4$ is evidently cannot be an
invariant 2-torus. The foliation of $W^c$ into closed invariant curves
homotopic to zero implies that near fixed point $\hat{O}$ we also get
closed invariant curves homotopic to zero. Embedding of $W^c$ to $\T^4$
can give an immersed plane or an immersed cylinder. We think the following
hypothesis is valid.
\begin{hypothesis}
In a transitive case, the image of $W^c$, $p(W^c)$, is a transitive surjection. This means that
$p$ is an inclusion (no two points whose image is the same) and the
closure of $p(W^c)$ coincides with $\T^4. $
\end{hypothesis}

\section{Classification of partially hyperbolic automorphisms}

At the study of partially hyperbolic symplectic automorphisms a natural question
arises: when two such automorphisms are topologically conjugate. Recall that two homeomorphisms
$f_1, f_2$ of a metric space $M$ are called topologically conjugate, if there exists
a homeomorphism $h:M \to M$ such that $h\circ f_1= f_2\circ h.$

Obviously, two automorphisms, one of which has transitive unstable foliation and another one
for which the unstable foliation generates the foliation $\mathbb{T}^4$ into two-dimensional tori
are topologically nonequivalent. However, two partially hyperbolic symplectic automor-phisms of
the torus having both transitive unstable foliations can also be nonequivalent.

Note, that classification of ergodic automorphisms of the torus from a measure theory point of
view is given by their entropy \cite{Katz}. This follows from the Ornstein isomorphism theorem
\cite{Orn} and the fact that every ergodic automorphism of a torus is Bernoulli one with
respect to the Lebesgue measure \cite{Katz}. But we want to obtain a
topological classification. The existence of a conjugating homeomorphism
$h: \T^4 \to \T^4$ for two automorphisms $f_A, f_{A'}$ implies the
relation $H\circ A = A'\circ H$ in the fundamental group $\Z^4$ of the
torus, here $H$ is the linear homomorphism in $\Z^4$ generated by $h$.
Thus, matrices $A,A'$ are similar by $H$. Matrix $H$ is integer with its
determinant $\pm 1.$ This is not occasional due to the following Arov's
theorem \cite{Arov}. We formulate the consequence of the Arov's theorem
\begin{theorem}
Two ergodic automorphisms $T,P$ of compact commutative group $X$ are
topolo-gically conjugate if and only if they are isomorphic, that is,
an isomorphism $Q:X\to X$ exists such that $Q\circ T=P\circ Q.$
\end{theorem}
Thus, the conjugacy problem for two transitive (i.e. ergodic) automorphisms
$f_A$, $f_{A'}$ is reduced to the problem: when two integer symplectic
matrices $A, A'$ are integer unimodularly conjugate. A necessary condition
is the equality of their characteristic polynomials with integer coefficients.
Therefore, $A$ and $A'$ should have the same companion matrix.
As was proved above, matrices $A$ and its companion matrix are
similar by means of some integer unimodular matrix. Thus, we get
similarity of $A, A'$ and hence conjugacy of $f_A$ and $f_{A'}.$
So, the following classification theorem for the transitive case is valid.
\begin{theorem}
	Let $A$ be a partially hyperbolic integer unimodular matrix generating a transitive
automorphism $f_A$ on $\T^4$. Then $f_A$ is symplectic w.r.t. a non-standard symplectic
structure. $f_A$ is integer unimodular conjugate to the automorphism $f_B$,	where $B$ is
the companion matrix of the characteristic polynomial of the matrix $A$.
\end{theorem}

It follows from above considerations the theorem on the structure of a
decomposable symplectic automorphism.
\begin{theorem}
If $f_A$ is a decomposable symplectic automorphism of $\T^4$, then $f_A$ is
conjugated by a diffeomorphism $f_S$ generated by an integer unimodular
matrix $S$ in $\R^4$ with the direct product of two automorphism $f_H$ and
$f_I$ on $\T^2\times \T^2$ given by a hyperbolic matrix $H$ and a
periodic matrix $I$ such that $I^k = E$ for $k\in\{3,4,6\}$.
\end{theorem}

\section{Acknowledgement}

The work by L.L. was supported by the Russian Science Foundation under the grant
19-11-00280. K.T. thanks for the partial financial support the Russian
Foundation for Basic Research under the grants 18-29-10081 and
19-01-00607.

\end{document}